\documentclass[a4paper,11pt]{article}
\usepackage{amsmath,amssymb,amsthm}
\theoremstyle{plain}
\newtheorem{thm}{Theorem}[section]

\newtheorem{cor}[thm]{Corollary}
\theoremstyle{definition}
\newtheorem{defn}[thm]{Definition}

\theoremstyle{remark}
\newtheorem{rem}{Remark}[thm]
\begin{document}
\title{SOME EXAMPLES OF SELF-SIMILAR 
SOLUTIONS AND TRANSLATING SOLITONS
 FOR LAGRANGIAN MEAN CURVATURE FLOW}
\author{Hiroshi Nakahara\\
                                    \\
Department of Mathematics \\
Tokyo Institute of Technology \\
2-21-1  O-okayama, Megro \\
Tokyo 152-8551 Japan \\
                                           \\
\textit{E-mail address}: nakahara.h.ab@m.titech.ac.jp}
\maketitle
\begin{abstract}
We construct new examples of self-similar solutions and translating solitons for Lagrangian mean
 curvature flow by extending the method of Joyce, Lee and Tsui \cite{MR2629511}. Those examples include examples in which the Lagrangian angle is arbitrarily
  small as the examples of Joyce, Lee and Tsui \cite{MR2629511}.
\end{abstract}
\section{Introduction}
\quad In recent years  the Lagrangian mean curvature flow has been 
extensively studied, as it is a key ingredient in the 
Strominger-Yau-Zaslow Conjecture and Thomas-Yau Conjecture.  Strominger-Yau-Zaslow Conjecture
 explains Mirror Symmetry of Calabi-Yau 3-folds. In Joyce, Lee and Tsui \cite{MR2629511}, many examples of self-similar 
 solutions and translating solitons for Lagrangian mean curvature flow are constructed. Those 
 Lagrangian submanifold $L$ are the total space of a 1-parameter family of quadrics $Q_s$, $s\in I$, where $I$ is an open interval 
 in $\mathbb{R}.$ 
 In this paper, we construct new examples of those Lagrangian submanifolds that link up with the examples of Lagrangian submanifolds given in \cite{MR2252892}, \cite{MR2015549}, \cite{MR2629511}, \cite{MR2563738} and so on. To do so we improve theorems in \cite{MR2629511} by describing Lagrangian submanifolds of the forms
 of Ansatz 3.1 and Ansatz 3.3 in \cite{MR2629511}.

 Let $L$ be a Lagrangian submanifold in $\mathbb{C}^n $. Define a function  $ c:L\to \mathbb{C}$ by the relation that $$\mathrm{d}z_1\wedge \cdots\wedge \mathrm{d}z_n |_L\equiv c\cdot \mathrm{vol}_L$$  where $\mathrm{vol}_L$ is the volume form of $L$.  Then $|c|\equiv 1$  holds. So  we can define 
{\it Lagrangian angle} $\theta :L \to \mathbb{R}$ or $\theta :L\to \mathbb{R}/2\pi \mathbb{Z}$  by the relation that $$\mathrm{d}z_1\wedge \cdots\wedge \mathrm{d}z_n |_L\equiv e^{i\theta}\mathrm{vol}_L.$$ On a Lagrangian submanifold $L$ in $\mathbb{C}^n$, the mean curvature vector $H$ is given by
\begin{equation}\label{H}
H=J\nabla \theta ,
\end{equation}
where $\nabla$  is the gradient on $L$ and $J$ is the standard complex structure in $\mathbb{C}^n$. The proof of \eqref{H} is given in \cite{MR1957663}.

\begin{defn}
Let $L\subset \mathbb{R}^N$ be a submanifold in $\mathbb{R}^N$.  $L$ is called a
\textit{self-similar solution} if  $H\equiv \alpha F^\perp $ on $L$ for some
 constant $\alpha \in \mathbb{R}$, where $F^\perp $ is the orthogonal projection of the position vector $F$ in $\mathbb{R}^N$ to 
 the normal bundle of $L$, and $H$ is the mean curvature vector of $L$ in $\mathbb{R}^N.$ It is called a \textit{self-shrinker} if $\alpha <0$ and a \textit{self-expander} if $c>0.$
 On the other hand $L\subset \mathbb{R}^N$ is called a \textit{translating soliton} if there exists a constant vector $T$ in $\mathbb{R}^N$ such that 
$H\equiv T^{\perp }$, where $T^\perp $ is the orthogonal projection of the constant vector $T$ in $\mathbb{R}^N$ to 
 the normal bundle of $L$ and $H$ is the mean curvature vector of $L$ in $\mathbb{R}^N.$ We call $T$ a {\it translating vector}.
\end{defn}
It is well-known that if $F$ is a self-similar solution then $F_t = \sqrt{2\alpha t}F $ is moved by the mean curvature flow, and 
if $F$ is a translating soliton then $F_t = F+tT $ is also moved by the mean curvature flow.

First we consider self-similar solutions.

\begin{thm}\label{main}
Let $\,B,\,C,\,\lambda_1,\cdots ,\lambda_n\in \mathbb{R}-\{0\},\, E>1,\quad a_1 ,\cdots ,a_n>0,$ and 
$ \alpha,\psi_1 ,\cdots ,\psi_n \in\mathbb{R}$
 be constants. Let $I\subset \mathbb{R}$ be a connected open neighborhood of $\,0\in \mathbb{R}$ such that 
$E \{\prod_{k=1}^{n}(1+2a_kB\lambda_k  s )\}  e^{2B\alpha s}-1$ and $1/a_j+ 2\lambda _j Bs$ are positive for any $j\in\{1,\cdots ,n\}$
 and any $s\in I.$ Define $r_1 ,\cdots ,r_n :I\to \mathbb{R} $ by 
 $$r_j =\sqrt{\frac{1}{a_j} + 2\lambda _j Bs}.$$ Define $\phi_1 ,\cdots ,\phi_n :I\to \mathbb{R} $ by 
\begin{equation*}
\phi _j=  \psi  _j +\int_{0}^{s}  \frac{\lambda _j |B|}{(\frac{1}{a_j} +2\lambda _j Bt)\sqrt{E\{\prod_{k=1}^{n}(1 + 2a_k B\lambda _k t )\}  e^{2B\alpha t}-1}} \mathrm{d}t.
\end{equation*}
Then the submanifold $L$ in $\mathbb{C}^n$  given by
\[ L=\{(x_1r_1(s) e^{i\phi_1(s)},\cdots,x_n r_n(s)e^{i\phi_n(s)})|\sum_{j=1}^n  \lambda_j x_j^2 =C, x_j \in\mathbb{R}, s\in I \} \]
is an immersed Lagrangian submanifold, and its position vector $F$ and mean curvature vector $H$ satisfy $CH\equiv \alpha  F^{\perp }$.
\end{thm}
\begin{rem}
In the situation of Theorem \ref{main}, let $I'\subset I$ be a subinterval such that $r_j$ has positive lower and upper bounds. Put $\alpha \neq 0, n\geq 2,\lambda_j >0 $ for $1\leq j\leq k<n$ and $\lambda_j <0$ for $k<j\leq n$ where $k$ is a positive integer less than $n.$ Let $t\in\mathbb{R}$ be a constant. 
Define $$L_t =\{(x_1r_1(s) e^{i\phi_1(s)},\cdots,x_n r_n(s)e^{i\phi_n(s)})|\sum_{j=1}^n  \lambda_j x_j^2 =2\alpha t, x_j \in\mathbb{R}, s\in I' \}.$$
The fact that the varifold $\cup_t \,L_t,\,-\infty<t<\infty, $ forms an eternal solution for Brakke flow without mass loss is proved similarly to Lee and Wang \cite{MR2563738}. So see Lee and Wang \cite{MR2563738}.
By Theorem \ref{main}, $L_t$ is a Lagrangian self-shrinker if $ t<0,$ a Lagrangian self-expander if $ t>0,$ and a Lagrangian cone if $t=0.$
\end{rem}
\begin{thm}\label{th3}
Let  $a_j>0$  , $\psi_j \in \mathbb{R}$,@$E>1$, and  $\alpha \geq 0$  be constants.
Define $r_j (s) :\mathbb{R}\to \mathbb{R}$  by  $r_j (s)= \sqrt{\frac{1}{a_j} +s^2}$.
Define  $\phi _j (s):\mathbb{R}\to \mathbb{R}$  by
\begin{equation}\label{int} 
\phi_j(s) =\psi _j + \int_{0}^{s} \frac{|t|}{(\frac{1}{a_j} +t^2)\sqrt{E \{\prod_{k=1}^{n}(1+a_k  t^2)\}  e^{\alpha t^2}-1}} \mathrm{d}t.
\end{equation}
Then the submanifold $L$ in $\mathbb{C}^n$  given by
\[ L=\{(x_1r_1(s) e^{i\phi_1(s)},\cdots,x_n r_n(s)e^{i\phi_n(s)})|\sum_{j=1}^n  x_j^2 =1, x_j \in\mathbb{R}, s\in\mathbb{R},\,s\neq 0\} \]
is an immersed Lagrangian diffeomorphic to $(\mathbb{R}-\{0\})\times S^{n-1},$ and its position vector $F$ and mean curvature vector $H$ satisfy $H\equiv \alpha F^{\perp }$.
\end{thm}
\begin{rem}
The manifold obtained as $\lim_{E\to 1+0} L$ is the same as Theorem C in \cite{MR2629511} . 
 So the condition $s\neq 0$ on the definition of $L$ is not necessary if $E=1$. 
If we put $E=1$ then changing $0\mapsto -\infty$ in the integral of \eqref{int} gives the example of \cite{MR2015549}.  
 \end{rem}
\begin{rem}\label{small}
Define $\bar{\phi}_1 ,\cdots ,\bar{\phi}_n >0$ by
$$\bar{\phi}_j =\int_{0}^{\infty} \frac{|t|}{(\frac{1}{a_j} +t^2)\sqrt{E \{\prod_{k=1}^{n}(1+a_k  t^2)\}  e^{\alpha t^2}-1}} \mathrm{d}t.$$ We put $\alpha >0.$ From the proof of Theorem \ref{th3} and \eqref{theta} in {\S}~\ref{S2} the Lagrangian angle $\theta$ satisfies 
\begin{equation}\label{section2}
\begin{split}
\theta        &=\sum_{j} \phi_j +\mathrm{arg}(s +i\frac{|s|}{\sqrt{E \{\prod_{k=1}^{n}(1+a_k  s^2)\}  e^{\alpha s^2}-1}}) \quad \mathrm{and}\\
\dot{\theta}&=\frac{-\alpha |s|}{\sqrt{E \{\prod_{k=1}^{n}(1+a_k  s^2)\}  e^{\alpha s^2}-1}} .
\end{split}
\end{equation}
It follows that  
$\theta$ is strictly decreasing. We define the submanifold $L_1 \subset L$ by restricting $s>0$ and $L_2 \subset L$ by restricting $s<0$. Therefore we have $L=L_1 \cup L_2 .$ We rewrite $\theta_1 , \theta_2$ as the Lagrangian angle of $L_1 , L_2$
 respectively. Then $\lim_{s\to +\infty}\theta_1(s)<\theta_1(s)<\lim_{s\to0+0}\theta_1(s)\quad \mathrm{and} \quad \lim_{s\to 0-0}\theta_2(s)<\theta_2(s)<\lim_{s\to -\infty}\theta_2(s).$ So from the first equation of \eqref{section2} we have 
$$\sum_{j} \psi_j +\sum_{j}\bar{\phi}_j <\theta_1 <\sum_{j} \psi_j + \mathrm{tan}^{-1}\frac{1}{\sqrt{E-1}}$$ and 
 $$\sum_{j} \psi_j +\pi - \mathrm{tan}^{-1}\frac{1}{\sqrt{E-1}}<\theta_2 <\sum_{j} \psi_j + \pi -\sum_{j}\bar{\phi}_j.$$
 Therefore by choosing $\mathrm{tan}^{-1}(1/\sqrt{E-1})$ close to $0,$ that is, choosing $E$ close to $\infty$, the oscillation of the Lagrangian angle of $L_1 ,L_2$ can be made arbitrarily small. 
 Furthermore $\Phi :(a_1 ,\cdots ,a_n)\mapsto (\bar{\phi}_1 ,\cdots ,\bar{\phi}_n)$ gives a diffeomorphism
\begin{multline*} 
\Phi :(0,\infty)^n \rightarrow \{(\bar{\phi}_1 ,\cdots ,\bar{\phi}_n)\in (0,\mathrm{tan}^{-1}\frac{1}{\sqrt{E-1}})^n|\, 0<\sum_{j}\bar{\phi}_j<\\
\mathrm{tan}^{-1}\frac{1}{\sqrt{E-1}}\}.
\end{multline*}
We can prove that $\Phi$ is a diffeomorphism similarly to the proof of Theorem D in \cite{MR2629511}. So we omit the proof. Therefore by choosing $\sum_{j}\bar{\phi}_j$ close to $\mathrm{tan}^{-1}(1/\sqrt{E-1}),$ the oscillation of the Lagrangian angle of $L_1 ,L_2$ can also be made arbitrarily small. 
\end{rem}
\begin{rem}
If we put $B=1/2,\,C=\lambda_1 =\cdots =\lambda_n =1,\,\alpha \geq 0$ in the situation of Theorem \ref{main} then $L$ is the same as $L$ in Theorem
 \ref{th3} where $s>0.$
\end{rem}
Next we turn to translating solitons.
\begin{thm}\label{main2}
Let $\,B,\,\lambda_1,\cdots ,\lambda_{n-1}\in \mathbb{R}-\{0\},\, E>1,\quad a_1 ,\cdots ,a_{n-1}>0,$ and 
$ \alpha,\,\psi_1 ,\cdots ,\psi_{n-1} \in\mathbb{R},\,K\in \mathbb{C}$
 be constants. Let $I\subset \mathbb{R}$ be a connected open neighborhood of $\,0\in \mathbb{R}$ such that 
$E \{\prod_{k=1}^{n-1}(1+2a_k B\lambda_k  s )\}  e^{2B\alpha s}-1$ and $1/a_j + 2\lambda _j Bs$ are positive for any $j\in\{1,\cdots ,n-1\}$
 and any $s\in I.$ Define $r_1 ,\cdots ,r_{n-1} :I\to \mathbb{R} $ by 
 $$r_j =\sqrt{\frac{1}{a_j} + 2\lambda _j Bs}.$$ Define $\phi_1 ,\cdots ,\phi_{n-1} :I\to \mathbb{R} $ by 
 $$\phi _j=  \psi  _j +\int_{0}^{s}  \frac{\lambda _j |B|}{(\frac{1}{a_j} +2\lambda _j Bt)\sqrt{E\{\prod_{k=1}^{n-1}(1 + 2a_k B\lambda _k t )\}  e^{2B\alpha t}-1}} \mathrm{d}t.$$
Then the submanifold $L$ in $\mathbb{C}^n$ given by
\[ L=\{(x_1 r_1(s)e^{i\phi_1 (s)},\cdots,x_{n-1}r_{n-1}(s)e^{i\phi_{n-1} (s)},-\frac{1}{2}\sum_{j=1}^{n-1}\lambda_j x_j^2 + Bs+\]
\[i|B|\int_0^s \frac{\mathrm{d}t}{\sqrt{E\{\prod_{k=1}^{n-1}(1 + 2a_k B\lambda _k t )\}  e^{2B\alpha t}-1}}+K)|
x_1,\cdots,x_{n-1}\in \mathbb{R},s\in I \}  \]
is an immersed Lagrangian submanifold, and its mean curvature vector $H$ satisfy $H\equiv T^{\perp }$ where $T= (0,\cdots ,0,\alpha)\in \mathbb{C}^n$.
\end{thm}
\begin{thm}\label{last}
Let  $a_1 ,\cdots ,a_{n-1}>0$  , $\psi_1 ,\cdots ,\psi_{n-1} \in \mathbb{R}$,@$E>1$, and  $\alpha \geq 0$  be constants.
Define $r_j (s) :\mathbb{R}\to \mathbb{R}$  by  $r_j (s)= \sqrt{\frac{1}{a_j} +s^2}$.
Define  $\phi _j (s):\mathbb{R}\to \mathbb{R}$  by
\[\phi_j(s) =\psi _j + \int_{0}^{s} \frac{|t|}{(\frac{1}{a_j} +t^2)\sqrt{E \{\prod_{k=1}^{n-1}(1+a_k  t^2)\}  e^{\alpha t^2}-1}} \mathrm{d}t.\]
Then the submanifold $L$ in $\mathbb{C}^n$  given by
\[ L=\{(x_1r_1(s) e^{i\phi_1(s)},\cdots,x_{n-1} r_{n-1}(s)e^{i\phi_{n-1}(s)},-\frac{1}{2}\sum_{j=1}^{n-1} x_j^2 + \frac{1}{2}s^2+\]
\[i\int_0^s \frac{|t|\mathrm{d}t}{\sqrt{E\{\prod_{k=1}^{n-1}(1 + a_k t^2 )\}  e^{\alpha t^2}-1}})| x_1 \cdots ,x_{n-1} \in\mathbb{R}, s\in\mathbb{R},\,s\neq 0\} \]
is an immersed Lagrangian diffeomorphic to $(\mathbb{R}-\{0\})\times \mathbb{R}^{n-1}$
 , and its mean curvature vector $H$ satisfy $H\equiv T^{\perp }$ where $T=(0,\cdots ,0,\alpha)\in \mathbb{C}^n$.
\end{thm}
\begin{rem}
If we put $\psi_1 =\cdots =\psi_{n-1}=0$ then the manifold obtained as $\lim_{E\to 1+0} L$ is the same as Theorem G in \cite{MR2629511} . 
 So the condition $s\neq 0$ on the definition of $L$ is not necessary if $E=1$.
 \end{rem}
\begin{rem}
We define the submanifold $L_1 \subset L$ by restricting $s>0$ and $L_2 \subset L$ by restricting $s<0$. Similarly to Remark \ref{small} if we fix $\alpha >0$ then 
by choosing the parameters $a_1 ,\cdots ,a_n >0$ or $E>1,$ the oscillation of the Lagrangian angle of $L_1 , L_2$ can be made arbitrarily small. 
\end{rem}
\begin{rem}
If we put $B=1/2,\,\lambda_1 =\cdots =\lambda_{n-1} =1,\,K=0,\, \alpha \geq 0$ in the situation of Theorem \ref{main2} then $L$ is the same as $L$ in Theorem
 \ref{last} where $s>0.$
\end{rem}
\section{Proofs for self-similar solutions}\label{S2}
In order to prove Theorem \ref{main} and Theorem \ref{th3} we use the following Theorem \ref{th1} and Theorem \ref{th2} that are a slight generalization of Theorem A and Theorem B in \cite{MR2629511}.
The following Theorem \ref{th1} sets up the ordinary differential equations
 for immersed Lagrangian submanifolds diffeorphic to $ \mathcal{S}^1 \times \mathcal{S}^{m-1} \times 
\mathbb{R}^{n-m}$ or $ \mathbb{R} \times \mathcal{S}^{m-1} \times 
\mathbb{R}^{n-m}$ where $1\leq m \leq n.$
\begin{thm}\label{th1}
Let $I\subset  \mathbb{R}$ be an open interval. Let $ \lambda_1,\cdots ,\lambda_n ,C\in \mathbb{R}-\{0\}$ be  constants, and $\omega_1 ,\cdots ,\omega_n:I\rightarrow \mathbb{C}-\{0\}$ and $f:I\to \mathbb{C}-\{0\}     $ be smooth functions. Suppose that
\begin{equation}\label{d}
\frac{\mathrm{d}\omega_j}{\mathrm{d}s}=\frac{\lambda_j f}{\overline{\omega_j}}     ,\quad j=1,\cdots,n       
\end{equation}
hold in $I$.
Then the submanifold $L$ in $\mathbb{C}^n$ given by
\begin{equation}\label{l}
L=\{(x_1\omega_1(s),\cdots,x_n\omega_n(s))|\sum_{j=1}^n \lambda_j x_j^2 =C, x_j \in\mathbb{R}, s \in I     \} 
\end{equation}
is a Lagrangian submanifold in $\mathbb{C}^n$, with Lagrangian angle 
\begin{equation}\label{theta}
\theta(s)=\arg (\omega_1\cdots \omega_n f)
\end{equation}
 at  $(x_1\omega_1(s),\cdots ,x_n\omega_n(s))\in L,$
so $\theta $ is a function depending only on $s$, not on $x_1,\cdots,x_n     .   $ 
Further we have 
\[H=\frac{\dot{\theta}}{\langle \frac{\partial}{\partial s}, \frac{\partial}{\partial s} \rangle     }J\left( \frac{\partial}{\partial s}  \right)    \]
and
\[
F^\bot   =-\frac{C\, \mathrm{Im}(f)}{ \langle \frac{\partial}{\partial s}, \frac{\partial}{\partial s} \rangle     }  J\left( \frac{\partial}{\partial s}  \right)                 .   \]
\end{thm}
\begin{rem}
A direct calculation shows that $$\langle \frac{\partial}{\partial s}, \frac{\partial}{\partial s} \rangle =\sum_{j=1}^{n} \frac{\lambda_j^2 x_j^2 |f|^2}{|\omega_j|^2 }.$$
\end{rem}
\begin{rem}
Let $\tilde{\theta}:I\to \mathbb{R} $ be a function satisfying $$\frac{\mathrm{d}\tilde{\theta}}{\mathrm{d}s}=-\alpha \, \mathrm{Im}(e^{i\tilde{\theta}}\,\overline{\omega_1 \cdots \omega_n}).$$
 Put $$f=e^{i\tilde{\theta}}\,\overline{\omega_1 \cdots \omega_n}.$$
Then it is proved in \cite{MR2629511} that $L$ is a Lagrangian submanifold and $L$ is a self-similar solution. This is Theorem A in \cite{MR2629511}.
\end{rem}
\begin{cor}\label{cor0913}
Let $I\subset  \mathbb{R}$ be an open interval. Let $ \lambda_1,\cdots ,\lambda_n ,C\in \mathbb{R}-\{0\}$ be  constants, and $\omega_1 ,\cdots ,\omega_n:I\rightarrow \mathbb{C}-\{0\}$ and $f:I\to \mathbb{C}-\{0\}     $ be smooth functions. Suppose that
\begin{equation*}
\frac{\mathrm{d}\omega_j}{\mathrm{d}s}=\frac{\lambda_j f}{\overline{\omega_j}}     ,\quad j=1,\cdots,n       
\end{equation*}
and
\begin{equation*}
\frac{\mathrm{d}}{\mathrm{d} s}\arg (\omega_1 \cdots \omega_n f) =-\alpha \mathrm{Im}(f)
\end{equation*}
hold in $I$.
Then the submanifold $L_C$ in $\mathbb{C}^n$ given by
\begin{equation*}
L_C=\{(x_1\omega_1(s),\cdots,x_n\omega_n(s))|\sum_{j=1}^n \lambda_j x_j^2 =C, x_j \in\mathbb{R}, s \in I     \} 
\end{equation*}
is a Lagrangian submanifold in $\mathbb{C}^n,$ and its position vector $F$ and mean curvature vector $H$ satisfy $CH\equiv \alpha F^{\perp }.$ Moreover, If $|\omega_j|$ has positive lower and upper bounds and $\alpha \in \mathbb{R}-\{0\}, n\geq 2,\lambda_j >0 $ for $1\leq j\leq k<n$ and $\lambda_j <0$ for $k<j\leq n$ where $k$ is a positive integer less than $n$  then the varifold $\cup_t \,L_{2\alpha t},\,-\infty<t<\infty, $ forms an eternal solution for Brakke flow without mass loss. The fact is proved similarly to Lee and Wang \cite{MR2563738}. Then $L_{2\alpha t}$ is a Lagrangian self-shrinker if \,$ t<0,$ a Lagrangian self-expander if \,$ t>0,$ and a Lagrangian cone if \,$t=0.$
\end{cor}

The following Theorem \ref{th2} gives the solution to ordinary deferential equations in Theorem \ref{th1}. 
\begin{thm}\label{th2}
In the situation of Theorem \ref{th1}  write     $\omega_j := r_j e^{i \phi_j}$, for functions $r_1,\cdots ,r_n:I\to (0,\infty), \phi_1 ,\cdots ,\phi_n:I\to \mathbb{R}$.
Fix $s_0 \in I.$ Define $u:I \to \mathbb{R}$   by
\[    u(s):=2\int_{s_0}^s \mathrm{Re}(f(t)) \mathrm{d}t               .                               \]
Then we have
$$r_j^2 =\alpha_j + \lambda_j u ,\quad\phi_j =\psi_j + \int_{s_0}^s \frac{\lambda_j\, \mathrm{Im}(f(t))}{\alpha_j + \lambda_j u(t)}\mathrm{d}t,  $$ 
with $\alpha_j =r_j^2(s_0), \quad    \psi  _j =\phi _j (s_0)             . $
\end{thm}
\begin{rem}
Thus if $f$ is explicitly given, the ordinary differential equation \eqref{d} is solved by Theorem \ref{th2}. 
\end{rem}

\textit{Proof of Theorem \ref{th1}.}
If we set $$x_n =\sqrt{\frac{1}{\lambda _n}(C-\lambda_1{x_1}^2-\cdots - \lambda_{n-1}{x_{n-1}}^2}),$$ then 
$(x_1,\cdots,x_{n-1},s)$ is a coordinate of $L.$ With this coordinate we have
\begin{equation*}
\begin{split}
    & (\mathrm{d}z_1\wedge \mathrm{d}\overline{z}_1   +\cdots + \mathrm{d}z_n\wedge \mathrm{d}\overline{z}_n    )|_L\\
=  &(\omega_1 \mathrm{d}x_1 + x_1 \dot{\omega}_1\mathrm{d}s)\wedge(\overline{\omega_1 }\mathrm{d}x_1 + x_1 \overline{\dot{\omega}_1}\mathrm{d}s)+\cdots +\\        
    & (\omega_{n-1} \mathrm{d}x_{n-1} + x_{n-1} \dot{\omega}_{n-1}\mathrm{d}s)\wedge(\overline{\omega_{n-1} }\mathrm{d}x_{n-1} + x_{n-1} \overline{\dot{\omega}_{n-1}}\mathrm{d}s)+\\
    &(-\frac{\lambda_1 x_1}{\lambda_n x_n}\omega_n\mathrm{d}x_1-\cdots -\frac{\lambda_{n-1}x_{n-1}}{\lambda_n x_n}\omega_n  \mathrm{d}x_{n-1}  +x_n\dot{\omega}_n\mathrm{d}s   )   \wedge \\         
    &(-\frac{\lambda_1 x_1}{\lambda_n x_n} \overline{\omega_n}\mathrm{d}x_1- \cdots -\frac{\lambda_{n-1}x_{n-1}}{\lambda_n x_n} \overline{\omega_n}  \mathrm{d}x_{n-1}  +x_n \overline{\dot{\omega}_n}\mathrm{d}s ) \\
=   & 2\mathrm{Im}(x_1\omega_1 \overline{\dot{\omega}_1})\mathrm{d}x_1 \wedge \mathrm{d}s
+\cdots +2\mathrm{Im}(x_{n-1}\omega_{n-1} \overline{\dot{\omega}_{n-1}})\mathrm{d}x_{n-1}\wedge \mathrm{d}s \\
    &  +2\mathrm{Im}(-\frac{\lambda_1}{\lambda_n}x_1 \omega_n \overline{\dot{\omega}_n})\mathrm{d}x_1 \wedge \mathrm{d}s+\cdots + 2\mathrm{Im}(-\frac{\lambda_{n-1}}{\lambda_n}x_{n-1} \omega_n \overline{\dot{\omega}_n})\mathrm{d}x_{n-1}\wedge \mathrm{d}s\\
=  &0,
\end{split}
\end{equation*}
and
\begin{equation*}
\begin{split}
    & \qquad \mathrm{d}z_1 \wedge \cdots \wedge \mathrm{d}z_n |_L\\
    & =
\begin{vmatrix}
\omega_1 & \cdots & \cdots  & x_1\dot{\omega}_1\\
             & \ddots &          0         & \vdots         \\
      0       &            &   \omega_{n-1}&x_{n-1}\dot{\omega}_{n-1} \\
-\frac{\lambda_1 x_1}{\lambda_n x_n}\omega_n  & \cdots  &  -\frac{\lambda_{n-1} x_{n-1}}{\lambda_n x_n}\omega_n  
& x_n \dot{\omega}_n                    &        
\end{vmatrix}
\begin{matrix}
\mathrm{d}x_1\wedge \cdots \wedge \mathrm{d}x_{n-1}\wedge \mathrm{d}s 
\end{matrix}\\
  &=\frac{\omega_1\cdots \omega _n}{\lambda _n x_n}(\sum_{j=1}^{n}\frac{\lambda _j x_j^2 \dot{\omega}_j}{\omega _j})\mathrm{d}x_1\wedge \cdots \wedge \mathrm{d}x_{n-1}\wedge \mathrm{d}s \\
   & =\frac{\omega_1\cdots \omega _nf}{\lambda _n x_n} \sum_j \frac{\lambda _j^2 x_j^2}{|\omega _j|^2}\mathrm{d}x_1\wedge \cdots \wedge \mathrm{d}x_{n-1}\wedge \mathrm{d}s.
\end{split}
\end{equation*}
It follows that $L$ is a nonsingular immersed Lagrangian, with Lagrangian angle 
$$\theta(s) =\arg (\omega_1\cdots \omega_n f) $$ at $(x_1\omega_1(s),\cdots ,x_n \omega_n (s)).$
Since 
$$\mathrm{d}\theta=\dot{\theta}\mathrm{d}s,$$   
and 
\begin{equation*}
\begin{split}
\langle \frac{\partial }{\partial s}, \frac{\partial }{\partial x_j}\rangle   
         &= \langle (x_1\dot{\omega}_1(s),\cdots,x_n \dot{\omega}_n(s)),(0,\cdots,0, \omega_j ,0,\cdots,0,- \frac{\lambda_j x_j}{\lambda_n x_n}\omega_n) \rangle \\
         &=\mathrm{Re}(x_j \dot{\omega}_j\overline{\omega_j} -\frac{\lambda_j}{\lambda_n}x_j \dot{\omega}_n \overline{\omega_n})\\
         &=x_j \lambda_j \,\mathrm{Re}(f-f)  \\
          &=0
\end{split}
\end{equation*}
for any $1\leq j\leq n-1,$ 
it follows that 
$$\nabla \theta =\frac{\dot{\theta}}{\langle \frac{\partial}{\partial s}, \frac{\partial}{\partial s} \rangle     }\left( \frac{\partial}{\partial s}  \right)  . $$
Therefore we obtain
$$H=\frac{\dot{\theta}}{\langle \frac{\partial}{\partial s}, \frac{\partial}{\partial s} \rangle     }J\left( \frac{\partial}{\partial s}  \right)   .$$
The normal projection of the position vector $F$ is computed by
\begin{equation*}
\begin{split}
\langle F,J \frac{\partial}{\partial x_j} \rangle 
    &= \langle (x_1 \omega_1,\cdots,x_n \omega_n),i(0,\cdots,0, \omega_j ,0,\cdots,0,- \frac{\lambda_j x_j}{\lambda_n x_n}\omega_n) \rangle \\
    &=0,
\end{split}
\end{equation*}
and
\begin{equation*}
\begin{split}
\langle F,J \frac{\partial}{\partial s} \rangle 
   & =\langle (x_1 \omega_1,\cdots,x_n \omega_n),i(x_1\dot{\omega}_1 (s),\cdots,x_n \dot{\omega}_n(s)) \rangle\\    
   & =\mathrm{Re}(\sum {x_j }^2\overline{\omega_j}\cdot i\dot{\omega}_j) \\
   &=-C\,\mathrm{Im}(f)      .
\end{split}
\end{equation*}
It follows that   $$  F^\bot   =-\frac{C\,\mathrm{Im}(f)}{ \langle \frac{\partial}{\partial s}, \frac{\partial}{\partial s} \rangle     }  J\left( \frac{\partial}{\partial s}  \right)   .$$
This completes the proof of Theorem \ref{th1}.\qed

\textit{Proof of Theorem \ref{th2}.}
Since we have
$$\frac{\mathrm{d} r_j^2}{\mathrm{d}s}=\frac{\mathrm{d}}{\mathrm{d}s} (\omega _j \overline{\omega _j})=\lambda _j f+ \lambda _j \bar{f}=2\lambda _j \,\mathrm{Re}(f),$$
we obtain $$r_j^2 =\alpha_j + \lambda_j u . $$ 
So we  have $$\dot{r}_j=\frac{\lambda _j \,\mathrm{Re}{f}}{r_j}.$$
By our assumption we have 
$$\dot{\omega}_j=\frac{\lambda_j f}{\overline{\omega_j}}.$$
This is equivalent to
 $$ \dot{r}_j   e^{i\phi _j }+ r_j i\dot{\phi}_j e^{i\phi _j }=\frac{\lambda _j f}{r_j e^{-i \phi _j}}  . $$ 
Then we have
$$\frac{\lambda _j \,\mathrm{Re}(f)}{r_j}+r_j i \dot{\phi}_j=\frac{\lambda _j f}{r_j}.$$
Then
$$ \dot{\phi}_j =\frac{\lambda _j \,\mathrm{Im}(f)}{\alpha _j +\lambda _j u} .$$
Therefore
$$\phi_j =\psi_j + \int_{s_0}^s \frac{\lambda_j \,\mathrm{Im}(f(t))}{\alpha_j + \lambda_j u(t)}\mathrm{d}t.$$
This completes the proof of Theorem \ref{th2}.
\qed

From Theorem \ref{th1} \,if \,$\dot{\theta}=-\alpha \,\mathrm{Im}(f)  $ then $L$ is a self-similar solution.
 From Theorem \ref{th2} and \eqref{theta}, $\dot{\theta}=-\alpha \,\mathrm{Im}(f)  $ is equivalent to
\begin{equation}\label{eq:1}    
\sum_j \frac{\lambda_j \,\mathrm{Im}(f)}{\alpha_j + 2\lambda_j \int_{s_0}^s \mathrm{Re}(f)\mathrm{d}t}  + \frac{\mathrm{d}}{\mathrm{d}s}\mathrm{arg}(f)= -\alpha \mathrm{Im}(f) .
\end{equation}
Therefore if \eqref{eq:1} holds then $L$ is a self-similar solution. So we can get self-similar solutions by getting solutions of \eqref{eq:1}. For example if we put $f\equiv i$ 
then $f$ is a solution of \eqref{eq:1}. In this case if we put $\alpha_1 =\cdots =\alpha _n =0, \,\lambda_1 ,\cdots ,\lambda_n \in\mathbb{Z}-\{0\}$ then the construction reduces to that of Lee and Wang \cite{MR2563738}.

\textit{Proof of Theorem \ref{main}.}
Define $\omega_j :I\to \mathbb{C}-\{0\}$ by $ \omega_j =r_j e^{i\phi_j}.$ Define 
$f :I\to \mathbb{C}-\{0\}$ by 
$$f= \frac{\overline{\omega_j}\,\dot{\omega}_j}{\lambda_j}.$$ A direct calculation shows that \[f=B+i \frac{|B|}{\sqrt{E\{\prod_{k=1}^{n}(1 + 2a_kB\lambda _k s)\}  e^{2B\alpha s}-1}}  .  \]
Then $f\neq 0 $ in $I.$ Apply Theorem 2.1 to the data $\omega_j , f $ above.
Then we get a Lagrangian submanifold $L$ defined by \eqref{l}. A direct calculation shows that $f$ satisfies \eqref{eq:1}. So this completes the proof. 
\qed

\textit{Proof of Theorem \ref{th3}.}
Put $I=\mathbb{R}_{>0}$ or $I=\mathbb{R}_{<0}.$ Define $\omega_j :I\to \mathbb{C}-\{0\}$ by $ \omega_j =r_j e^{i\phi_j}.$ Define 
$f :I\to \mathbb{C}-\{0\}$ by 
$$f= \overline{\omega_j}\,\dot{\omega}_j.$$ A direct calculation shows that $$f=s+i\frac{|s|}{\sqrt{E \{\prod_{k=1}^{n}(1+a_k  s^2)\}  e^{\alpha s^2}-1}}.$$
Then $f\neq 0 $ in $I.$ Apply Theorem 2.1 to the data $\omega_j , f $ above and $\lambda_1 =\cdots =\lambda_n =C=1.$ 
Then we get a Lagrangian submanifold $L$ defined by \eqref{l}. A direct calculation shows that $f$ satisfies \eqref{eq:1}. So this completes the proof. 
\qed

\section{Proofs for translating solitons}
In order to prove Theorem \ref{main2} and Theorem \ref{last} we use the following Theorem \ref{th4} and Corollary \ref{cor} that are a slight generalization of Theorem A and Theorem B in \cite{MR2629511}.
The following Theorem \ref{th4} sets up the ordinary differential equations
 for immersed Lagrangian submanifolds diffeorphic to $ \mathbb{R}^n .$
\begin{thm}\label{th4}
Fix $n\geq 2$. Let  $\lambda_1 ,\cdots , \lambda_{n-1} \in \mathbb{R} \backslash \{0\}$ and  $\alpha  \in  \mathbb{R}$  be constants, $I$ be an open 
interval in $\mathbb{R}$, and $\omega_1,\cdots ,\omega_{n-1}:I\to \mathbb{C}\backslash \{0\} $ and $\beta :I\to \mathbb{C}$
 be smooth functions. Suppose that
\begin{equation}\label{d2} 
\frac{\mathrm{d}\omega _j}{\mathrm{d}s}=\frac{\lambda_j }{\overline{\omega_j}}\cdot\frac{\mathrm{d}\beta}{\mathrm{d}s}
\end{equation}
 and
 \[\frac{\mathrm{d}\beta}{\mathrm{d}s}\neq0 \]
hold in $I$. Then the submanifold $L$ in $\mathbb{C}^n$ given by
\begin{equation}\label{Lt}
\begin{split}
 L=\{ (x_1 \omega_1(s),\cdots,x_{n-1}&\omega_{n-1}(s),-\frac{1}{2}\sum_{j=1}^{n-1}\lambda_j x_j^2 + \beta (s))|\\
                                                    &x_1,\cdots,x_{n-1}\in \mathbb{R},s\in I \}  
\end{split}
\end{equation}
is an immersed Lagrangian submanifolds diffeomorphic to $\mathbb{R}^n$, with Lagrangian angle 
\begin{equation}\label{theta2}
\theta (s)= \mathrm{arg} (\omega_1 \cdots
 \omega_{n-1} \dot{\beta} ).
\end{equation}
Further we have 
\[ H= \frac{\dot{\theta}}{\langle \frac{\partial}{\partial s}, \frac{\partial}{\partial s} \rangle}J\left( \frac{\partial}{\partial s}  \right)  \] and
\[ T^{\perp }=\frac{-\alpha \,\mathrm{Im}(\dot{\beta})}{\langle \frac{\partial}{\partial s}, \frac{\partial}{\partial s} \rangle   }J\left( \frac{\partial}{\partial s}  \right)     , \]
where $T=(0,\cdots ,0,\alpha)\in \mathbb{C}^n $ is a constant vector.
\end{thm}
\begin{rem}
A direct calculation shows that $$\langle \frac{\partial}{\partial s}, \frac{\partial}{\partial s} \rangle =\sum_{j=1}^{n-1} \frac{\lambda_j^2 x_j^2 |\dot{\beta}|^2}{|\omega_j|^2 }+|-\frac{1}{2}\sum_{j=1}^{n-1}\lambda_j x_j^2 +\beta|^2 .$$
\end{rem}
\begin{rem}
Let $\tilde{\theta}:I\to \mathbb{R} $ be a function satisfying $$\frac{\mathrm{d}\tilde{\theta}}{\mathrm{d}s}=-\alpha \, \mathrm{Im}(e^{i\tilde{\theta}}\,\overline{\omega_1 \cdots \omega_{n-1}}).$$
Suppose that $$\frac{\mathrm{d}\beta}{\mathrm{d}s}=e^{i\tilde{\theta}}\,\overline{\omega_1 \cdots \omega_{n-1}}.$$ hold in $I.$ 
Then it is proved in \cite{MR2629511} that $L$ is a Lagrangian submanifold and $L$ is a translating soliton with translating vector $(0,\cdots ,0,\alpha)\in \mathbb{C}^n $. This is Theorem G in \cite{MR2629511}.
\end{rem}
The following Corollary \ref{cor} gives the solution to ordinary deferential equations in Theorem \ref{th4}. 
\begin{cor}\label{cor}
In the situation of Theorem \ref{th4},
write   $\omega_j := r_j e^{i \phi_j}$, for functions $r_1,\cdots ,r_{n-1}:I\to (0,\infty), \phi_1 ,\cdots ,\phi_{n-1}:I\to \mathbb{R}$.
Define $u:I \to \mathbb{R}$   by
\[    u(s):=2\int_{s_0}^s \mathrm{Re}(\dot{\beta}(t)) \mathrm{d}t    .       \]
Then we have 
$$r_j^2 =\alpha_j + \lambda_j u, \quad \phi_j =\psi_j + \int_{s_0}^s \frac{\lambda_j \,\mathrm{Im}(\dot{\beta}(t))}{\alpha_j + \lambda_j u(t)}\mathrm{d}t  $$with $\alpha_j =r_j^2(s_0), \quad \psi  _j =\phi _j (s_0)      .  $
\end{cor}
\begin{rem}
Thus if $\beta$ is explicitly given, the ordinary differential equation \eqref{d2} is solved by Theorem \ref{cor}. 
\end{rem}

\textit{Proof of Theorem \ref{th4}.}
We consider the map
\[\iota :\mathbb{R}^{n-1}\times I\to \mathbb{C}^n \]
 given by
\[\iota((x_1,\cdots ,x_{n-1}),s)=(x_1 \omega_1(s),\cdots,x_{n-1}\omega_{n-1}(s),\frac{1}{2}\sum_{j=1}^{n-1}\lambda_j x_j^2 + \beta (s)).\]
Then we have $L=\iota (\mathbb{R}^{n-1}\times I).$ 
So $(x_1,\cdots,x_{n-1},s)$ is a coordinate of $L.$ With this coordinate we have
\[
\begin{split}
   &  (\mathrm{d}z_1 \wedge \mathrm{d}\overline{z}_1+\cdots +\mathrm{d}z_n\wedge \mathrm{d}\overline{z}_n)|_L  \\
= &(\omega_1 \mathrm{d}x_1 +x_1 \dot{\omega}_1 \mathrm{d}s)\wedge (\overline{\omega_1} \mathrm{d}x_1 +x_1 
\overline{\dot{\omega}_1}\mathrm{d}s)+\cdots +\\
   &(\omega_{n-1} \mathrm{d}x_{n-1} +
\dot{\omega}_{n-1}\mathrm{d}s)\wedge (\overline{\omega_{n-1}} \mathrm{d}x_{n-1}+ x_{n-1}
\overline{\dot{\omega}_{n-1}}\mathrm{d}s)\\
    &+(-\lambda_1 x_1 \mathrm{d}x_1-\cdots
-\lambda_{n-1}x_{n-1}\mathrm{d}x_{n-1}+\dot{\beta}(s)\mathrm{d}s)\\
    &\wedge (-\lambda_1 x_1 \mathrm{d}x_1-\cdots
-\lambda_{n-1}x_{n-1}\mathrm{d}x_{n-1}+\overline{\dot{\beta}(s)}\mathrm{d}s)\\
=   &2i\mathrm{Im}(x_1 \omega_1 \overline{\dot{\omega}_1})\mathrm{d}x_1 \wedge
\mathrm{d}s+\cdots +2i\mathrm{Im}(x_{n-1} \omega_{n-1} \overline{\dot{\omega}_{n-1}})\mathrm{d}x_{n-1} \wedge \mathrm{d}s\\
    &-2i\mathrm{Im}(\lambda_1 x_1 \overline{\dot{\beta}})\mathrm{d}x_1 \wedge \mathrm{d}s-\cdots 
-2i\mathrm{Im}(\lambda_{n-1} x_{n-1} \overline{\dot{\beta}})\mathrm{d}x_{n-1} \wedge \mathrm{d}s\\
=   &0
\end{split}
\]
and
\[
\begin{split}
    &  \mathrm{d}z_1 \wedge \cdots \wedge \mathrm{d}z_n|_L   \\
=  & (\omega_1 \mathrm{d}x_1+x_1 \dot{\omega}_1\mathrm{d}s)\wedge \cdots \wedge (\omega_{n-1} \mathrm{d}x_{n-1}+x_{n-1} \dot{\omega}_{n-1}\mathrm{d}s)\wedge \\
    &(-\lambda_1 x_1 \mathrm{d}x_1 -\cdots -\lambda_{n-1} x_{n-1}\mathrm{d}x_{n-1} +\dot{\beta}\mathrm{d}s) \\
=   &\omega_1 \cdots \omega_{n-1}\dot{\beta}(1+ \sum_{j=1}^{n-1}\frac{\lambda_j^2 x_j^2}{|\omega_j|^2})\mathrm{d}x_1 \wedge \cdots \wedge \mathrm{d}x_{n-1}\wedge \mathrm{d}s .
\end{split}
\]
It follows that $L$ is a nonsingular immersed Lagrangian, with Lagrangian angle 
\[\theta =\mathrm{arg}(\omega_1 \cdots \omega_{n-1}\dot{\beta}).\]
Since 
\[ 
\begin{split}
    & \langle \frac{\partial }{\partial s}, \frac{\partial }{\partial x_j}\rangle \\
=  &\langle (x_1\dot{\omega}_1(s),\cdots, x_{n-1}\dot{\omega}_{n-1}(s),\dot{\beta}),(0,\cdots,0, \omega_j ,0,\cdots,0,- \lambda_j x_j ) \rangle \\
=  &\mathrm{Re}(\overline{\omega_j}x_j \dot{\omega}_j -\lambda_j x_j \dot{\beta})\\
=  &0,
\end{split}
\]
and
$\mathrm{d}\theta = \dot{\theta}\mathrm{d}s,$ 
it follows that $$\nabla \theta =\frac{\dot{\theta}}{\langle \frac{\partial}{\partial s}, \frac{\partial}{\partial s} \rangle     }\left( \frac{\partial}{\partial s}  \right)  . $$ Therefore we obtain $$H=\frac{\dot{\theta}}{\langle \frac{\partial}{\partial s}, \frac{\partial}{\partial s} \rangle     }J\left( \frac{\partial}{\partial s}  \right).   $$
The normal projection of the position vector $F$ is computed by
\[
\begin{split}
   & \langle T,J \frac{\partial}{\partial x_j} \rangle \\
 = &\langle (0,\cdots,0,\alpha),
i(0,\cdots,0, \omega_j ,0,\cdots,0,-\lambda_j x_j) \rangle \\
=  &0
\end{split}
\]
and
\[
\begin{split}
\langle T,J \frac{\partial}{\partial s} \rangle =  &\langle (0,\cdots,0,\alpha),i(x_1\dot{\omega}_1(s),\cdots,x_{n-1} \dot{\omega}_{n-1}(s),\dot{\beta}) \rangle \\
                                                              =  &\mathrm{Re}(\alpha i \dot{\beta})\\
                                                              =  &-\alpha \, \mathrm{Im}(\dot{\beta})     .
\end{split}
\]
It follows that   $$  T^\bot   =-\frac{\alpha \,\mathrm{Im}(\dot{\beta})}{ \langle \frac{\partial}{\partial s}, \frac{\partial}{\partial s} \rangle     }  J\left( \frac{\partial}{\partial s}  \right)   .$$
This completes the proof of Theorem \ref{th4}.\qed

As the proof of Theorem G in \cite{MR2629511} there is an another proof of Theorem \ref{th4} which obtains 
Theorem \ref{th4} from Theorem \ref{th1} by a limiting procedure. See Joyce, Lee and tsui \cite{MR2629511}.

\textit{Proof of Corollary \ref{cor}.}
Changing $f$ to $\dot{\beta}$ in the proof of Theorem \ref{th2} gives the proof of Corollary \ref{cor}.
\qed

From Theorem \ref{th4} \,if \,$\dot{\theta}=-\alpha \,\mathrm{Im}(\dot{\beta})  $ then $L$ is a translating soliton with translating vector 
$T=(0,\cdots ,0,\alpha)\in \mathbb{C}^n $.
 From Theorem \ref{cor} and \eqref{theta2}, $\dot{\theta}=-\alpha \,\mathrm{Im}(\dot{\beta})  $ is equivalent to
\begin{equation}\label{eq:4}   
\sum_{j=1}^{n-1} \frac{\lambda_j \,\mathrm{Im}(\dot{\beta})}{\alpha_j + 2\lambda_j \int_{s_0}^s \mathrm{Re}(\dot{\beta})\mathrm{d}t}  + \frac{\mathrm{d}}{\mathrm{d}s}\mathrm{arg}(\dot{\beta})= -\alpha \,\mathrm{Im}(\dot{\beta}) .
\end{equation}
Therefore if \eqref{eq:4} holds then $L$ is a translating soliton with translating vector 
$T=(0,\cdots ,0,\alpha)\in \mathbb{C}^n $.

\textit{Proof of Theorem \ref{main2}.}
Define $\omega_j :I\to \mathbb{C}-\{0\}$ by $ \omega_j =r_j e^{i\phi_j}.$ Define 
$\beta:I\to \mathbb{C}-\{0\}$ by 
$$\dot{\beta}= \frac{\overline{\omega_j}\,\dot{\omega}_j}{\lambda_j}\,\,\,\mathrm{and} \,\,\,\beta(0)=-K.$$ A direct calculation shows that \[\beta (s)=\int_0^s(B+i \frac{|B|}{\sqrt{E\{\prod_{k=1}^{n-1}(1 + 2a_kB\lambda _k s)\}  e^{2B\alpha s}-1}} )\mathrm{d}s +K .  \]
Then $\dot{\beta}\neq 0 $ in $I.$ Apply Theorem \ref{th4} to the data $\omega_j , \beta $ above.
Then we get a Lagrangian submanifold $L$ defined by \eqref{Lt}. A direct calculation shows that $\beta$ satisfies \eqref{eq:4}. So this completes the proof. 
\qed

\textit{Proof of Theorem \ref{last}.}
Put $I=\mathbb{R}_{>0}$ or $I=\mathbb{R}_{<0}.$ Define $\omega_j :I\to \mathbb{C}-\{0\}$ by $ \omega_j =r_j e^{i\phi_j}.$ Define 
$\beta :I\to \mathbb{C}-\{0\}$ by 
$$\dot{\beta}= \overline{\omega_j}\,\dot{\omega}_j\,\,\,\mathrm{and}\,\,\, \beta (0)=0.$$ A direct calculation shows that $$\beta (s)=\int_0^s ( s+i\frac{|s|}{\sqrt{E \{\prod_{k=1}^{n-1}(1+a_k  s^2) \}  e^{\alpha s^2}-1}}\,)\mathrm{d}s.$$
Then $\dot{\beta}\neq 0 $ in $I.$ Apply Theorem \ref{th4} to the data $\omega_j , \beta $ above and $\lambda_1 =\cdots =\lambda_{n-1} =1.$ 
Then we get a Lagrangian submanifold $L$ defined by \eqref{Lt}. A direct calculation shows that $\beta$ satisfies \eqref{eq:4}. So this completes the proof. 
\qed

\subsection*{Acknowledgements}
The author would like to thank the supervisor Akito Futaki. He also wishes to thank Masataka Shibata, Yuji Terashima, Mitutaka 
Murayama and Kota Hattori for useful conversations.

\bibliographystyle{plain}
\bibliography{ref2}
\end{document}